\documentclass[notitlepage,leqno,10pt]{article}

\usepackage{latexsym}
\usepackage{amsmath}
\usepackage{amsfonts}
\usepackage{amssymb}
\usepackage{color}

\renewcommand{\a }{\alpha }

\renewcommand{\d}{\delta }
\newcommand{\D }{\Delta }

\newcommand{\n }{\nabla }

\newcommand{\intbar}{\mathop{\int\makebox(-13.5,0){\rule[4pt]{.7em}{0.3pt}}%
\kern-6pt}\nolimits}

\newcommand{\be}{\begin{equation}}
\newcommand{\ee}{\end{equation}}

\newenvironment{pf}{\noindent{\sc Proof}.\enspace}{\rule{2mm}{2mm}\medskip}
\newenvironment{pfn}{\noindent{\sc Proof}}{\rule{2mm}{2mm}\medskip}

\newcommand{\R}{\mathbb{R}}

\def\bM{\partial M}
\def\dn{\partial_{\nu}}

\def\gt{\tilde{g}}

\def\sdtu{\sigma_2(\gt^{-1} A^1_{\gt})}
\def\sdu{\sigma_2(g^{-1} A^1_{g})}
\def\sutu{\sigma_1(\gt^{-1} A^1_{\gt})}
\def\suu{\sigma_1(g^{-1} A^1_{g})}
\def\sigut{\tilde{\sigma_1}}
\def\sigu{\sigma_1}
\def\sigdt{\tilde{\sigma_2}}
\def\sigd{\sigma_2}

\def\be{\begin{equation}}
\def\ee{\end{equation}}
\def\bea{\begin{eqnarray*}}
\def\eea{\end{eqnarray*}}
\def\f{\frac}

\author{Giovanni CATINO$^{a}$ and Cheikh Birahim NDIAYE$^{b}$}

\date{}

\title{\bf Integral pinching results for  manifolds with boundary}

\begin{document}

\parindent=0pt

\newtheorem{lem}{Lemma}[section]
\newtheorem{pro}[lem]{Proposition}
\newtheorem{thm}[lem]{Theorem}
\newtheorem{rem}[lem]{Remark}
\newtheorem{cor}[lem]{Corollary}
\newtheorem{df}[lem]{Definition}
\newtheorem{claim}[lem]{Claim}
\newtheorem{conj}[lem]{Conjecture}

\maketitle

\begin{center}

\

{\small

\noindent $^a$ SISSA -- International School for Advanced Studies

Via Beirut 2--4,

I-34014 Trieste - Italy
}

\

{\small

\noindent $^b$ Mathematisches Institut der Universit\"at T\"ubingen

Auf der Morgenstelle 10,

D-72076 T\"ubingen - Germany

}

\end{center}

\footnotetext[1]{E-mail addresses: catino@sissa.it, ndiaye@everest.mathematik.uni-tuebingen.de}

\

\

\begin{center}
{\bf Abstract}

\end{center}
We prove that some Riemannian manifolds with boundary satisfying an explicit integral pinching condition are spherical space forms. More precisely, we show that three-dimensional Riemannian manifolds with totally geodesic boundary, positive scalar curvature and an explicit integral pinching between the $L^2$-norm of the scalar curvature and the $L^2$-norm of the Ricci tensor are spherical space forms with totally geodesic boundary. Moreover, we prove also that four-dimensional Riemannian manifolds with umbilic boundary, positive Yamabe invariant and an explicit integral pinching between the total integral of the $(Q,T)$-curvature and the $L^2$-norm of the Weyl curvature are spherical space forms with totally geodesic boundary. As a consequence, we show that a certain conformally invariant operator, which plays an important role in Conformal Geometry, is non-negative and has trivial kernel if the Yamabe invariant is positive and verifies a pinching condition together with the total integral of the $(Q,T)$-curvature. As an application of the latter spectral analysis, we show the existence of conformal metrics with constant $Q$-curvature, constant $T$-curvature, and zero mean curvature  under the latter assumptions.
\begin{center}

\bigskip\bigskip


\bigskip

\centerline{\bf AMS subject classification:   53C24, 53C20,
53C21, 53C25}

\end{center}

\centerline{}

\section{Introduction}\label{s:intro}

One of the most important questions about the relation between algebraic properties of the full curvature tensor
and the topology of manifolds is under which conditions on its curvature tensor a Riemannian manifold is homeomorphic or diffeomorphic to a space of constant sectional curvature, namely a space form. A model example is the classical {\em sphere theorem}  conjectured by Rauch \cite{ra}, and  which says that any closed, simply connected and $\frac{1}{4}$-pinched Riemannian manifold is diffeomorphic to the standard sphere. The topological version was proved by Berger \cite{be} and Klingenberg \cite{kl}. Just recently the original conjecture has been settled by Brendle and Schoen \cite{bs}, using a result of Bohm and Wilking \cite{bw}.

On the other hand, many sphere like theorems appeared in the literature in the last $30$ years in connection to the celebrated Ricci flow. Just to mention some of them which are related to our results, we start by recalling the pioneering work of R. Hamilton \cite{ha}. Using the Ricci flow, he proved the following theorem
\begin{thm} ({\bf Hamilton})\\
If \;$(M, g)$ is a closed three--dimensional Riemannian manifold
with positive Ricci curvature, then $M$ is diffeomorphic to a spherical space form, i.e. $M$ admits a metric with constant positive sectional curvature.
\end{thm}

Later  C. Margerin $\cite{ma}$ proved an optimal curvature characterization of the smooth $4$-sphere. We recall Margerin' s theorem in a form where the optimality issue is not apparent, but enough for the link with our work. We define the weak pinching quantity $$
W P_g = \frac{|W_g|^ 2_g + 2 |E_g |^2_g}{ R^2_ g},
$$ 
where $W_g$ denoting  the Weyl tensor, $E_g$ the trace--free Ricci tensor and $| \cdot |_g$  the usual norm of a tensor with respect to the metric $g$. Here is the result
\begin{thm} ({\bf Margerin})\\
Let \;$(M, g)$ be a closed four--dimensional Riemannian manifold with positive scalar curvature. If the pinching condition $W P_g < \tfrac{1}{6}$ is satisfied, then $M$ is diffeomorphic to a spherical space form. Moreover, we get that the manifold $M$ is diffeomorphic to $S^4$ or $\R P^ 4$.
\end{thm}

\medskip

Much later, Chang, Gursky and Yang \cite{cgy2} proved a remarkable improvement of Margerin's theorem with assumptions which are in integral form, and conformally invariant too.
\begin{thm} ({\bf Chang-Gursky-Yang}) \\
 Let \;$(M, g)$ be a closed four--dimensional Riemannian manifold with positive Yamabe invariant. If the curvatures satisfy
$$
\int_M \left( |W_g |^2_g + 2 |E_g|_g^2 -\frac{1}{6}R^2_g\right) dV_g < 0,
$$
then $M$ is diffeomorphic to a spherical space form. Moreover, we get that the manifold $M$ is diffeomorphic to $S^4$ or $\R P^ 4$.
\end{thm}
Notice that the integral pinching condition can be written in the following form (for the definition of $Q_g$, see below)
$$
\int_M Q_gdV_g>\frac{1}{8}\int_M|W_g|^2 dV_g.
$$

\medskip

Recently, the first author and Z. Djadli \cite{cd} proved an integral pinching theorem in dimension three.
\begin{thm} ({\bf Catino-Djadli})\\
Let \;$(M, g)$ be a closed three--dimensional Riemannian manifold with positive scalar curvature. If
$$
\int_M |Ric_g|^2_g dV\leq \frac{3}{8}\int_M R^2_g dV_g ,
$$
then $M$ is diffeomorphic to a spherical space form.
\end{thm}

\medskip

On the other hand, the Ricci flow techniques have also been used to get sphere like theorems for manifolds with boundary. An example which is of interest to us is the following result of  Shen \cite{she}.
\begin{thm}\label{eq:shen} ({\bf Shen})\\
If \;$(M, g)$ is a compact three--dimensional Riemannian manifold with totally geodesic boundary and positive Ricci curvature, then $M$ admits a metric with constant positive sectional curvatures in the interior and totally geodesic boundary.
\end{thm}

Using the Ricci flow for manifolds with boundary defined by Shen \cite{she}, a very easy adaptation of the arguments of Margerin \cite{ma}, immediately yields the following theorem
\begin{thm}\label{eq:bmar}
Let $(M, g)$ be a compact four--dimensional Riemannian manifold  with totally geodesic boundary and  positive scalar curvature. If the pinching condition $W P_g < \tfrac{1}{6}$ is satisfied, then $M$ admits a metric with constant positive sectional curvatures in the interior and totally geodesic boundary.
\end{thm}

\medskip

Our goal in this paper is to provide counterparts of the results of Chang-Gursky-Yang and Catino-Djadli for manifolds with boundary. The first result we will prove is the following
\begin{thm}\label{T3_gen}
Let $(M,g)$ be a compact three--dimensional Riemannian
manifold with totally geodesic boundary and positive scalar curvature. If
$$
\int_M |Ric_{g}|^{2}_{g}\,dV_{g} \leq\f{3}{8}\int_M
R^{2}_{g}\,dV_{g}\,,
$$
then $M$ admits a metric with constant positive sectional curvatures in the interior and totally geodesic boundary.
\end{thm}

\medskip

In order to state our second result on four manifolds with boundary, we need to recall some notions from Conformal Geometry. We start by recalling the Paneitz operator and its associated curvature invariant called $Q$-curvature. In 1983, Paneitz has discovered a  conformally covariant differential operator on four dimensional compact smooth Riemannian manifolds with smooth boundary $(M,g)$ (see \cite{p1}). To this operator, Branson \cite{bo} has associated a natural curvature invariant called $Q$-curvature. They are defined in terms of Ricci tensor \;$Ric_{g}$\; and scalar curvature \;$R_{g}$\; of the manifold \;$(M,g)$\; as follows
\begin{equation*}
P^4_g\varphi=\D_{g}^{2}\varphi+div_{g}\left((\frac{2}{3}R_{g}g-2Ric_{g})d\varphi\right),\;\;\;\;\;\;\;Q_g=-\frac{1}{12}(\D_{g}R_{g}-R_{g}^{2}+3|Ric_{g}|^{2}),
\end{equation*}
where\;$\varphi$\; is any smooth function on \;$M$, \;$div_g$\;is the divergence and \;$d$\;is the De Rham differential.

\medskip

Similarly, Chang and Qing \cite{cq1}, have  discovered a boundary operator \;$P^3_g$\;defined on the boundary of compact four dimensional smooth Riemannian manifolds and a natural third-order curvature \;$T_g$\;associated to \;$P^3_g$ as follows
\begin{equation*}
P^3_g\varphi=\frac{1}{2}\frac{\partial {\D_g\varphi}}{\partial n_g}+\D_{\hat g}\frac{\partial \varphi}{\partial n_g}-2H_g\D_{\hat g}\varphi+(L_g)_{ab}(\nabla_{\hat g}\varphi)_{a}(\nabla_{\hat g}\varphi)_{b}+\nabla_{\hat g}H_g.\nabla_{\hat g}\varphi+(F-\frac{R_g}{3})\frac{\partial \varphi}{\partial n_g}.
\end{equation*}
\begin{equation*}
T_g=-\frac{1}{12}\frac{\partial R_g}{\partial n_g}+\frac{1}{2}R_gH_g-<G_g,L_g>+3H_g^3-\frac{1}{3}Tr(L^3)+\D_{\hat g}H_g,
\end{equation*}
where\;$\varphi$\; is any smooth function on \;$M$,\;\;\;$\hat g$\;is the metric induced by \;$g$\;on\;\;$\partial M$, \;$L_g=(L_g)_{ab}=-\frac{1}{2}\frac{\partial g_{ab}}{\partial n_g}$\\is the second fundamental form of \;$\partial M$,\;\;$H_g=\frac{1}{3}tr(L_g)=\frac{1}{3}g^{ab}L_{ab}$\;(\;$g^{a,b}$\;are the entries of the inverse \;$g^{-1}$\;of the metric\; $g$)\;is the mean curvature of \;$\partial M$,\;$R^k_{bcd}$\; is the Riemann curvature tensor\;\; $F=R^{a}_{nan}$,\;\;$R_{abcd}=g_{ak}R^{k}_{bcd}$\;(\;$g_{a,k}$\;are the entries of the metric \;$g$)\;and\;\;$<G_g,L_g>=R_{anbn}(L_g)_{ab}$, $\frac{\partial }{\partial n_g}$ is the inward normal derivative with respect to $g$. We recall that $(M,g)$ has umbilic boundary if $L_{g}=\lambda g$ for some constant $\lambda$. If $L_{g}=0$ we say that the boundary is totally geodesic.

\medskip

A remarkable property of the couple of operators $(P^4_g,P^3_g)$ is that, as  the couple Laplace-Beltrami operator and  Neumann operator governs the transformation law of the Gauss curvature and the geodesic curvature on compact surfaces with boundary under conformal change of metric,\;$(P^4_g,P^3_g)$\;does the same for \;$(Q_g,T_g)$\;on compact four dimensional smooth Riemannian manifolds with boundary. In fact, after a conformal change of metric \;$ g_u=e^{2u}g$\;we have that
\begin{equation}\label{eq:lawqt}
\left\{
     \begin{split}
P^4_{g_u}=e^{-4u}P^4_g;\\
P^3_{g_u}=e^{-3u}P^3_{g};
     \end{split}
   \right.
\qquad \mbox{and}\qquad
\left\{
\begin{split}
P^4_g u+2Q_g=2Q_{ g_u}e^{4u}\;\;\text{in }\;\;M\\
P^3_g u+T_g=T_{ g_u}e^{3u}\;\;\text{on}\;\;\partial M.
\end{split}
\right.
\end{equation}
An other very important role played by the couple of curvatures $(Q_g,T_g)$ in Conformal Geometry is that they arise in the well-known Gauss-Bonnet-Chern formula. More precisely
\begin{equation}\label{eq:gbc}
\int_{M}(Q_{g}+\frac{|W_{g}|^{2}}{8})dV_{g}+\oint_{\bM}(T+Z)dS_g=4\pi^{2}\chi(M)
\end{equation}
where \;$W_g$\; and \;$ZdS_g$\;(for the definition of \;$Z$\; see \cite{cq1}) are pointwise conformally invariant. Moreover, it turns out that \;$Z$\; vanishes when the boundary is totally geodesic. Setting
\begin{equation*}
\kappa_{P^4_g}=\int_{M}Q_gdV_g,\;\;\;\;\;\kappa_{P^3_g}=\oint_{\bM}T_gdS_{g},
\end{equation*}
from \;$\eqref{eq:gbc}$, thanks to the fact that \;$W_gdV_g$\;and \;$ZdS_g$\;are pointwise conformally invariant, we have that \;$\kappa_{P^4_g}+\kappa_{P^3_g}$\;is conformally invariant, and will be denoted by
\begin{equation}\label{eq:invc}
\kappa_{(P^4,P^3)}=\kappa_{P^4_g}+\kappa_{P^3_g}.
\end{equation}
In addition to the conformally invariant quantity $\kappa_{(P^4,P^3)}$ of a  compact four-dimensional Riemannian manifold with boundary,  there exists also the Yamabe invariant of the conformal class $[g]=\{\tilde g=e^{2u}g,\;\;u\in C^{\infty}(M)\}$ defined by the 
\begin{equation}\label{eq:yamabe}
Y(M,\partial M, [g])=\inf_{\tilde g\in[g], vol_{\tilde g}=1}\int_{M}R_{\tilde g}dV_{\tilde g}+\oint_{\partial M}H_{\tilde g}dS_{\tilde g}
\end{equation}
We recall that this invariant is defined for every compact Riemannian manifold with boundary of dimension greater or equal to $3$. 

\medskip

Now we are ready to state our result on four manifolds with boundary.
\begin{thm}\label{T4_gen} Let \;$(M,g)$ be a compact four--dimensional Riemannian
manifold with umbilic boundary. If \;$Y(M,\bM,[g])>0$ and if \;$\kappa_{(P^4,P^3)}>\frac{1}{8}\int_M |W_g|^2dV_g$, then \;$M$ admits a metric with constant positive sectional curvatures in the interior and totally geodesic boundary.\end{thm}

\medskip

The couple $(P^4_g,P^3_g)$ gives rise to an operator defined on $H_{\frac{\partial}{\partial n}}=\Big\{u\in H^2(M):\;\;\;\frac{\partial u}{\partial n_g}=0\Big\}$  whose spectral property is very important for uniformization problems on four manifolds with boundary. The latter operator that we denote by $P^{4,3}_g$ is defined as follows
\begin{equation*}
\left<P^{4,3}_gu,v\right>_{L^2(M)}=\int_{M}\left(\D_g u\D_gv+\frac{2}{3}R_g\nabla_g u\nabla_g v\right)dV_g-2\int_{M}Ric_g(\nabla_g u,\nabla_g v)dV_g-2\oint_{\bM}L_g( \nabla_{\hat g} u, \nabla_{\hat g} v)dS_g ,
\end{equation*}
for every \;$u,v\in H_{\frac{\partial }{\partial n}}$.

\medskip

As a byproduct of our analysis, we obtain the following spectral property for $P^{4,3}_g$.
\begin{thm}\label{eq:spect}
Let \;$(M,g)$ be a compact four--dimensional Riemannian manifold with umbilic boundary. Assuming \;$Y(M,\bM,[g])>0$ and \;$\kappa_{(P^4,P^3)}+\frac{1}{6}Y(M, \partial M, [g])^2>0$, then  \;$P^{4,3}_g$ is non--negative and \;$ker P^{4,3}_g\simeq\R$ .
\end{thm}

A direct consequence of Theorem $\ref{eq:spect}$ is the existence of constant $Q$-curvature and constant $T$-curvature conformal metrics on four--manifolds which verify the assumptions of Theorem $\ref{eq:spect}$
\begin{cor}\label{eq:existence}
Let \;$(M,g)$ be a compact four--dimensional Riemannian manifold with umbilic boundary. Assuming \;$Y(M,\bM,[g])>0$ and \;$\kappa_{(P^4,P^3)}+\frac{1}{6}Y(M, \partial M, [g])^2>0$, then \;$M$ carries a metric conformal to $g$ with constant \;$Q$-curvature, constant \;$T$-curvature and zero mean curvature.
\end{cor}

\

Proofs of Theorem \ref{T3_gen} and Theorem \ref{T4_gen} rely on the solution of some boundary value problems for fully nonlinear equations. Following \cite{gv} we will use the continuity method proving a priori estimates on the solutions to our equations. As a consequence of our work in dimension four, analising the spectral property of a certain operator, we will show that this operator is non--negative and with trivial kernel (Theorem \ref{eq:spect}). As a byproduct, we will prove then existence of conformal metrics with constant $Q$-curvature, constant $T$-curvature and zero mean curvature under certain conformally invariant assumptions (Corollary \ref{eq:existence}).

\

The plan of the paper is the following: in Section 2 we will introduce some notations, set up the boundary value problem; in Section 3 and 4 we will prove Theorem \ref{T3_gen} and Theorem \ref{T4_gen} on three and four manifolds respectively; finally Section 5 will be devote to the proof of Theorem \ref{eq:spect} and Corollary \ref{eq:existence}.

\

\section{Preliminaries and notations}
In this section, we give some notations and preliminaries like the notion of $k$-th symmetric elementary functions and some of their properties, the notion of $\sigma_k$-curvature of a Riemannian manifold, and some Moser-Trudinger type inequalities.\\
For this end, let $(M,g)$ be a compact, smooth, $n$--dimensional Riemannian manifold with boundary.
We will denote by $\nu_g$ the inner normal vector field with respect to the metric $g$ and by $\dn=\f{\partial}{\partial n_g}$ the inward normal derivative. Moreover $L_g$ and $H_g$ will be the second fundamental form $$L_{g,ab}=-\f{1}{2}\f{\partial g_{ab}}{\partial n_g},$$ and the mean curvature normalized, i.e. $$H_g=\f{1}{n-1}\,g^{ab}L_{g,ab}.$$

\

Given a section $A$ of the bundle of symmetric 2--tensors, we can use the metric to
raise an index and view $A$ as a tensor of type $(1,1)$, or
equivalently as a section of $End(TM)$. This allows us to define
$\sigma_k(g^{-1}A)$ the $k$-th elementary function of the
eigenvalues of $g^{-1}A$. More precisely we define

\begin{df} Let
$(\lambda_1,\cdots,\lambda_n)\in {\mathbb R}^n$. We view the
$k$-th elementary symmetric function as a function on ${\mathbb R}^n$:
$$
\sigma_{k}(\lambda_1,\cdots,\lambda_n)=\sum_{1\le i_1<\cdots<i_{k} \le
  n}\lambda_{i_1}\cdots\lambda_{i_{k}}\,,
$$
and we define
$$
\Gamma^+_{k}=\bigcap_{1\leq j\leq k}\{\sigma_j(\lambda_1,\cdots,\lambda_n)>0\}\subset {\mathbb
R}^n\,,
$$
\end{df}
For a symmetric linear transformation $A:V\rightarrow V$, where
$V$ is an $n$--dimensional inner product space, the notation
$A\in\Gamma^+_k$ will mean that the eigenvalues of $A$ lie in the
corresponding set. We note that this notation also makes sense for
a symmetric 2--tensor on a Riemannian manifold. If
$A\in\Gamma^+_k$, let
$\sigma_k^{1/k}(A)=\{\sigma_k(A)\}^{1/k}$.
\begin{df}\label{defNT} Let $A:V\rightarrow V$, where $V$ is an $n$--dimensional inner product
space. The $(k-1)$-th Newton transformation associated with
$A$ is
$$
T_{(k-1)}(A)=\sum_{j=0}^{k-1}(-1)^{k-1-j}\sigma_j(A)A^{k-1-j}\,.
$$
Also, for $t\in {\mathbb R}$ we define the linear transformation
$$
L^t(A)=T_{(k-1)}(A)+\f{1-t}{n-2}\sigma_1(T_{(k-1)}(A))\cdot I\,.
$$
\end{df}

We have the following list of properties (the proofs can be found in \cite{cns})
\begin{lem}\label{proper}\label{PropGamma}
\begin{enumerate}
\item[{\em{(i)}}] $\Gamma^+_k$ is an open
  convex cone with vertex at the origin, and we have the following
  sequence of inclusions
  $$
\Gamma_n^+\subset\Gamma_{n-1}^+\subset\cdots\subset\Gamma_1^+\,.
$$
\item[{\em{(ii)}}] If $A\in\Gamma^+_k$, then $T_{k-1}(A)$ is positive definite.
Hence for all $t\leq 1$, $L^t(A)$ is positive definite.
\item[{\em{(iii)}}] We have the identities
$$
T_{k-1}(A)^{ij} A_{ij}=k\,\sigma_k(A)\,,
$$
$$
T_{k-1}(A)^{ll}=(n-k+1)\sigma_{k-1}(A)\,.
$$
\item[{\em{(iv)}}] If $A\in\Gamma^+_k$, then
\bea
\sigma_{k-1}(A)\geq\f{k}{n-k+1}{n\choose k}^{1\over
  k}\sigma_k(A)^{(k-1)\over k}\,.
\eea
\item[{\em{(v)}}] If $A$ and $B$ are symmetric linear transformations,
$A,B\in\Gamma^+_k$, then $\forall\rho\in[0,1]$,
$\rho A+(1-\rho)B\in\Gamma^+_k,$ and
$$
\sigma_k^{\f{1}{k}}(\rho
A+(1-\rho)B)\geq\rho\sigma_k^{\f{1}{k}}(A)+(1-\rho)\sigma_k^{\f{1}{k}}(B)\,.
$$
In particular this gives the concavity of the function
$\sigma_k^{1\over k}$ in the cone $\Gamma_k^+$.
\end{enumerate}
\end{lem}

Next we give a Lemma about the variation of the $\sigma_k$ functional.
\begin{lem}\label{derform}
If $A \, : \, {\mathbb R}\rightarrow Hom(V,V)$, then
$$
\f{d}{ds}\sigma_k(A)(s)=\sum_{i,j}T_{(k-1)}(A)_{ij}(s)\f{d}{ds}(A)_{ij}(s)\,,
$$
i.e, the $(k-1)$-th Newton transformation is what arises when we
differentiate $\sigma_k$.
\end{lem}

We choose the tensor (here $t$ is a real number)
$$
A^t_g=\f{1}{n-2}\left(Ric_g-\f{t}{2(n-1)}R_g g\right)\,,
$$
where $Ric_g$ and $R_g$ denote the Ricci and the scalar curvature of
$g$ respectively. Note that for $t=1$, $A^1_g$ is the classical
Schouten tensor, namely $A^1_g=A_g:=\f{1}{n-2}\left(Ric_g-\f{1}{2(n-1)}R_g g\right)$, (see\cite {bes}). Hence,
with our notations,
$\sigma_k(g^{-1} A^t_g)$ denotes the $k$-th elementary symmetric
function of the eigenvalues of $g^{-1}A^t_g$.\\

Now, we give a Lemma which shows that metrics $g_1$, such that $A^t_{g_1}$ belong to the positive cone of order $2$, verify also additional pointwise algebraic inequalities. More precisely, we have
\begin{lem}\label{GVineq}
If for some metric $g_1$ on $M$ we have
$A_{g_1}^t\in\Gamma^+_2$, then
\bea
-A_{g_1}^t+\sigma_1({g_1^{-1}}A_{g_1}^t)g_1&>&0,\\
A_{g_1}^t+\f{n-2}{n}\sigma_1({g_1^{-1}}A_{g_1}^t)g_1&>&0.
\eea
\end{lem}

\

We will be concerned with the following equation for a conformal
metric $\gt=e^{-2u}g$:
\be\label{eqgen}
\left\{
\begin{array}{ll}
\sigma^{1/k}_k(g^{-1}A^t_{u})=fe^{2u}\hspace{1cm} & \mbox{in $M$},\\
\dn u=0\hspace{1cm} &\mbox{on $\bM$}.
\end{array}\right.
\ee
where $f$ is a positive function on $M$. Let
$\sigma_1(g^{-1}A^1_g)$ be the trace of $A^1_g$ with respect to
the metric $g$. We have the following formula for the
transformation of $A_g^t$ under this conformal change of metric:
\be\label{trasfz}
A^t_{\gt}=A^t_g+\nabla_g^2
u+\f{1-t}{n-2}(\Delta_g u) g+du\otimes du-\f{2-t}{2}|\nabla_g
u|_g^2 g\,.
\ee
Since
$$
A^t_g=A^1_g+\f{1-t}{n-2}\sigma_1(g^{-1}A^1_g)g\,,
$$
this formula follows easily from the standard formula for the
transformation of the Schouten tensor (\cite{via}):
\be\label{trasfz2} A^1_{\gt}=A^1_g+\nabla_g^2 u+du\otimes
du-\f{1}{2}|\nabla_g u|_g^2 g\,. \ee Using this formula we may write
~\eqref{eqgen} with respect to the background metric $g$
$$
\sigma_k\left(g^{-1}\left(A^t_g +\nabla_g^2 u+\f{1-t}{n-2}(\Delta_g
u)g+du\otimes du-\f{2-t}{2}|\nabla_g u|_g^2
g\right)\right)^{1/k}=f(x)e^{2u}\,.
$$

\

Now, we discuss the ellipticity properties of equation~\eqref{eqgen}.
\begin{pro}[Ellipticity property]\label{Inver}
Let $u\in C^2(M)$ be a solution of equation~\eqref{eqgen} for some
$t\leq 1$ and let ${\tilde g}=e^{-2u}g$. Assume that
$A^t_{\gt}\in\Gamma^+_k$. Then the linearized operator at
$u$, \\$\mathcal{L}^t:C^{2,\alpha}(M)\cap\{\dn u=0\,\,\mbox{on}\,\,\bM\}\rightarrow C^{\alpha}(M)$, is
elliptic and invertible $(0<\alpha<1)$.
\end{pro}
\begin{pf} Define the operator
$$
F_t[u, \nabla_g u, \nabla_g^2
u]=\sigma_k(g^{-1}A^t_{\gt})-f(x)^ke^{2ku}\,,
$$
so that solutions of the
equation~\eqref{eqgen} are exactly the zeroes of $F_t$. Define the
function $u_s=u+s\varphi$, then the linearization at $u$ of the
operator $F_t$ is defined by
\bea
\mathcal{L}^t(\varphi)&=&\f{d}{ds}F_t[u_s,\nabla_g u_s, \nabla_g^2
u_s]\Big|_{s=0}\\&=&\f{d}{ds}\left(\sigma_k(g^{-1}A^t_{\gt})\right)\Big|_{s=0}-\f{d}{ds}\left(f(x)^k
  e^{2ku_s}\right))\Big|_{s=0}\,.
\eea
>From Lemma~\ref{derform} we have
$$
\f{d}{ds}\left(\sigma_k(g^{-1}A^t_{\gt})\right))\Big|_{s=0}=T_{k-1}(g^{-1}A^t_{\gt})_{ij}\f{d}{ds}\left((A^t_{\gt}\right)_{ij}))\Big|_{s=0}\,.
$$
We compute
$$
\f{d}{ds}\left((A^t_{\gt}\right)_{ij}))\Big|_{s=0}=(\nabla^2_g\varphi)_{ij}+\f{1-t}{n-2}(\Delta_g
\varphi)g_{ij}-(2-t)\nabla_g u\cdot\nabla_g \varphi\,
g_{ij}+2du\otimes d\varphi\,.
$$
Easily we have also
$$
\f{d}{ds}\left(f(x)^k
  e^{2ku_s}\right))\Big|_{s=0}=2kf(x)^ke^{2ku}\,\varphi\,.
$$
Putting all together, we conclude
$$
\mathcal{L}^t(\varphi)=T_{k-1}(g^{-1}A^t_{\gt})_{ij}\left((\nabla^2_g\varphi)_{ij}+\f{1-t}{n-2}(\Delta_g
  \varphi)g_{ij}\right)-2kf(x)^ke^{2ku}\,\varphi+\cdots
$$
where the last terms denote additional ones witch are linear in
$\nabla_g\varphi$. The first term of the linearization is exactly the
one defined in~\ref{defNT},
i.e.
$$
L^t(A^t_{\gt})_{ij}=T_{k-1}(A^t_{\gt})_{ij}+\f{1-t}{n-2}T_{k-1}(A^t_{\gt})_{pp}\,\delta_{ij}\,.
$$
So finally, we have
$$\mathcal{L}^t(\varphi)=L^t(A^t_{\gt})_{ij}(\nabla^2_g\varphi)_{ij}-2kf(x)^ke^{2ku}\,\varphi+\cdots$$
Since $A^t_{\gt}\in\Gamma^+_k$, by Lemma~\ref{proper}, we have that
the tensor $L^t(A^t_{\gt})$ is positive definite. So, the linearized
operator at any solution $u$ must be elliptic. Note also that, by the
previous formula, the operator is of the form
$$
\mathcal{L}^t(\varphi)=E(\varphi)-c(x)\varphi\,,
$$ where $E(\varphi)$ is a second order linear elliptic operator and $c(x)$ is a strictly positive function on $M$, since $c(x)=2kf(x)^k e^{2ku}$ and $f(x)>0$. This allows us to invert this operator between the H\"older spaces $C^{2,\alpha}(M)\cap\{\dn u=0\,\,\mbox{on}\,\,\bM\}$ and $C^{\alpha}(M)$ (see for instance \cite{gt}).
\end{pf}

\

Next, we recall some Moser-Trudinger type inequalities which will be used to prove Corollary $\ref{eq:existence}$.
\begin{pro}\label{eq:mos-tru}
Assume  $(M,g)$ is a compact four-dimensional Riemannian manifold with boundary such that \;$P^{4,3}_{g}$\;is a non-negative operator with \;$Ker P^{4,3}_{g}\simeq\R$. Then  we have that for all \;$\alpha <16\pi^2$\;there exists a constant \;$C=C(M,g,\alpha)$\;such that
\begin{equation*}
\int_{M}e^{\frac{\alpha(u-\bar u)^2}{\left<P^{4,3}_{g}u,u\right>_{L^2(M)}}}dV_{g}\leq C,
\end{equation*}
for all\;$u\in H_{\frac{\partial }{\partial n}}$,\;and hence
\begin{equation*}
\log\int_{M}e^{4(u-\bar u)}\leq C+\frac{4}{\alpha}\left<P^{4,3}_{g}u,u\right>_{L^2(M)}\;\;\forall u\in H_{\frac{\partial }{\partial n}},
\end{equation*}
\end{pro}
where $\bar u=\frac{1}{Vol_g(M)}\int_M udV_g$, and $Vol_g(M)=\int_MdV_g$.\\

The latter Proposition can be found in \cite{nd1} together with its proof. The second inequality that we are going to state is a trace analogue of the previous one. Its proof can be found \cite{nd2}.
\begin{pro}\label{eq:mos-trub}
Assume \;$P^{4,3}_{g}$\;is a non-negative operator with \;$Ker P^{4,3}_{g}\simeq\R$. Then  we have that for all \;$\alpha <12\pi^2$\;there exists a constant \;$C=C(M,g,\alpha)$\;such that
\begin{equation}\label{eq:mts}
\oint_{\partial M}e^{\frac{\alpha(u-\bar u_{\partial M})^2}{\left<P^{4,3}_{g}u,u\right>_{L^2(M,g)}}}dS_{g}\leq C,
\end{equation}
for all\;$u\in H_{\frac{\partial }{\partial n}}$,\;and hence
\begin{equation}\label{eq:mts}
\log\oint_{\partial M}e^{3(u-\bar u_{\partial M})}dS_{g}\leq C+\frac{9}{4\alpha}\left<P^{4,3}_{g}u,u\right>_{L^2(M,g)}\;\;\forall u\in H_{\frac{\partial }{\partial n}}.
\end{equation}
\end{pro}
where $\bar u_{\partial M}=\frac{1}{Vol_{g}(\partial M)}\oint_{\partial M}udS_g$, and $Vol_g(\partial M)=\oint_{\partial M}dS_g$.\\

Now, we give a Lemma  (whose proof can be found in \cite{nd1}) which will be used  together with the above Moser-Trudinger type inequalities in order to prove Corollary $\ref{eq:existence}$. It says that under the assumptions \;$Ker P^{4,3}_g\simeq\R$\;and\;$P^{4,3}_g$\;non-negative, the map
\begin{equation*}
u\in H_{\frac{\partial }{\partial n}}\longrightarrow ||u||_{P^{4,3}_{g}}
=\left<P^{4,3}_{g}u,u\right>_{L^2(M)}^{\frac{1}{2}}
\end{equation*}
induces an equivalent norm to the standard norm of \;$H^{2}(M)$\;\;on \;\;$\{u\in H_{\frac{\partial }{\partial n}}\;\;\bar u=0\}$. More precisely we have the following
\begin{lem}\label{eq:neq}
Suppose \;$Ker P^{4,3}_g\simeq\R$\;and\;$P^{4,3}_g$\;non-negative then we have that \;$||\cdot||_{P^{4,3}_{g}}$\;is an equivalent norm to\;$||\cdot||_{H^2}$\;on  \;$\{u\in H_{\frac{\partial }{\partial n}}\;\;\bar u=0\}$
\end{lem}

\

Now we give a technical Lemma which will be used to prove the above theorems.
\begin{lem}\label{eq:boun}
Let $(M,g)$ be a compact $n$-dimensional Riemannian manifold with totally geodesic boundary. Assuming $u\in C^{2}(M)$ with $\frac{\partial u}{\partial n_g}=0$, then
$$
\frac{\partial |\n_g u|^2_g}{\partial n_g}=0,
$$
and
$$
A_g(\nu, \n_g u)=0.
$$
\end{lem}
\begin{pf}
First of all, using the fact that $\frac{\partial u}{\partial n_g}=0$, we derive
$$
|\n_gu|^2=g^{ab}\partial_{a}u\partial_{b}u.
$$
Thus we infer
$$
\frac{\partial (|\n_gu|^2)}{\partial n_g}=\frac{\partial g^{ab}}{\partial n_g}\partial_{a}u\partial_{b}u+2g^{ab}\frac{\partial(\partial_{a}u)}{\partial n_g}\partial_{b}u.
$$
Next, using the fact that $L_g=0$, one has $\frac{\partial g^{ab}}{\partial n_g}=0$.
 Moreover from the trivial identity $\frac{\partial(\partial_{a}u)}{\partial n_g}=\partial_{a}\left(\frac{\partial u}{\partial n_g}\right)$, we infer
$$
\frac{\partial(\partial_{a}u)}{\partial n_g}=0.
$$
Thus, we obtain
$$
\frac{\partial (|\n_g u|^2)}{\partial n_g}=0.
$$
This prove the first point. For the second one, we have
$$
A_g(\nu, \n_g u)=\frac{1}{n-2}\left(Ric_g(\nu,\n_g u)-\frac{1}{2(n-1)}R_g\frac{\partial u}{\partial n_g}\right).
$$
Thus, we get
$$
A_g(\nu, \n_g u)=\frac{1}{n-2}\,Ric_{\nu,a}\partial_{a}u.
$$
Now using the Codazzi-Mainardi equation, we get
$$
Ric_{\nu, a}=\n_{b}L_{g,ab}-\n_{a}H_g=0.
$$
So, we obtain.
$$
A_g(\nu, \n_g u)=0.
$$
This completes the proof of the lemma.
\end{pf}

\

\section{Three manifolds with boundary}

In this section, we present the proof of Theorem $\ref{T3_gen}$. We will prove a more general theorem so that Theorem $\ref{T3_gen}$ will be a direct corollary. In fact, we have
\begin{thm}\label{T3_1} Let $(M,g)$ be a compact three--dimensional Riemannian
manifold with totally geodesic boundary and positive scalar curvature.\\
There exists a positive constant $C=C(\mbox{diam}(M,g),\Vert \nabla^2 Rm\Vert)$ such that if
$$
\int_M\sigd(g^{-1}A^1_g)\,dV_g+C\left({7\over 10} -
  t_0\right)Y(M,[g])^2 >0\,,
$$
for some $t_0\leq 2/3$, then there exists a conformal metric
$\gt=e^{-2u}g$ with $R_{\gt}>0$,
$\sigma_2(g^{-1}A^{t_0}_{\gt})>0$ pointwise and totally geodesic boundary. Moreover, we have the
inequalities
\be\label{inequa}
(3t_0-2)R_{\gt}\gt<6Ric_{\gt}<3(2-t_0)R_{\gt}\gt\,.
\ee
\end{thm}

\

Throughout the sequel, $(M,g)$ will be a compact
3-dimensional Riemannian manifold with totally geodesic boundary and with positive scalar curvature.
Since $M$ is compact and $R_g>0$, there exists $t_0>\delta>-\infty$ such that $A^\delta_g$ is
positive definite (i.e. $Ric_g - {\delta \over 4}R_g g >0$ on $M$).
Note that $\delta$ only depends on $\Vert Rm_g\Vert$.

\

For $t\in[\delta, t_0]$, consider the path of
equations (in the sequel we use the notation
$A^t_{u_t} := A^t_{g_t}$ for $g_t$ given by
$g_t=e^{-2u_t}g$)
\be\label{eq3}
\left\{
\begin{array}{ll}
\sigma^{1/2}_2(g^{-1}A^t_{u_t})=fe^{2u_t}\hspace{1cm} & \mbox{in $M$},\\
\dn u_t=0\hspace{1cm} &\mbox{on $\bM$}.
\end{array}\right.
\ee
where
$f=\sigma^{1/2}_2(g^{-1}A^\delta_g)>0.$
Note that $u\equiv 0$ is a solution for $t=\delta$.

\

We use the continuity method. Define
$$\mathcal{S}=\left\{t\in[\delta, t_0]\mid\exists\mbox{\,a
solution\,}u_t\in
C^{2,\alpha}(M)\mbox{\,of\,}(\ref{eq3})\mbox{\,with\,}A^t_{u_t}\in\Gamma^+_2\right\}.$$
Clearly, with our choice of $f$, $u\equiv 0$ is a solution for
$t=\delta$. Since $A^{\delta}_g$ is positive definite, then
$\delta\in\mathcal{S}$. Hence $\mathcal{S}\neq\emptyset.$ Let
$t\in\mathcal{S}$, and $u_t$ be a solution. By Proposition
\ref{Inver}, the linearized operator at $u_t$,
$\mathcal{L}^t:C^{2,\alpha}(M)\cap\{\dn u=0\,\,\mbox{on}\,\,\bM\}\rightarrow C^{\alpha}(M)$, is
invertible. The implicit function theorem  tells us
that $\mathcal{S}$ is open. To prove that $\mathcal{S}$ is close we need to establish a priori $C^{2,\a}$ estimates for solutions of the equation (\ref{eq3}). To do this, we start by proving an upper bound estimate for solutions of $\eqref{eq3}$.

\

\begin{pro}[Upper bound]\label{UB3} Let $u_t\in C^2(M)$ be a solution of
(\ref{eq3}) for some $t\in[\delta, t_0]$. \\If $g_{t}=e^{-2u_{t}}g\in\Gamma_{2}^{+}$, then $u_t\leq\bar{\delta}$,
where $\bar{\delta}$ depends only on $\Vert Rm_g\Vert$.
\end{pro}
\begin{pf} From Lemma \ref{PropGamma} (iv), we have
$\sqrt{3}\sigma_2^{1/2}\leq\sigma_1$, so for all $p\in M$
$$\sqrt{3}fe^{2u_t}\leq\sigma_1(g^{-1}A^t_{u_t}).$$
Let $p\in M$ be a maximum of $u_t$. Since
the gradient terms vanish at $p$ (this is true also if $p\in\bM$, since $\dn u_t=0$ on $\bM$) we have
$(\Delta u_t)(p)\leq 0$. Then, using (\ref{trasfz}), we have
 \bea
\sqrt{3}f(p)e^{2u_t(p)}&\leq&\sigma_1(g^{-1}A^t_{u_t})(p)\\&=&\sigma_1(g^{-1}A^t_{g})(p)+(4-3t)(\Delta
u_t)(p)\\&\leq&\sigma_1(g^{-1}A^t_{g})(p)\\&\leq&\sigma_1(g^{-1}A^{\delta}_{g})(p).
\eea
Since $M$ is compact,
we have $u_t\leq\bar{\delta},$ for some $\bar{\delta}$ depending only on $\Vert Rm_g\Vert$.
\end{pf}

\

Next, we are going to show that solutions of $\eqref{eq3}$ which verify upper-bound estimates enjoy also gradient ones

\begin{pro}[Gradient estimate]\label{GE} Let $u_t\in C^3(M)$ be a solution
of (\ref{eq3}) for some $\delta\leq t\leq t_0$. Assume that
$u_t\leq\bar{\delta}$. Then $\parallel\nabla_g
u\parallel_{{g,\infty}}<C_1$, where $C_1$ depends only on
$\Vert\nabla Rm_g\Vert$ and $\bar\delta$.
\end{pro}
\begin{pf} Let $H:=|\nabla_g u|^2_g$. If the maximum of $H$ is in the interior, then $\nabla_g H=0$ and $\nabla^2_g H$ is negative semi-definite. If the maximum of $H$ is at the boundary, then by Lemma $\ref{eq:boun}$, $\frac{\partial H}{\partial n_g}=0$. Thus, we also have that $\nabla_g H=0$ and $\nabla^2_g H$ is negative semi-definite. Interior gradient estimates for equation (\ref{eq3}) were proved in \cite{gv} (Proposition 4.1). We remark that the same proof works for boundary gradient estimates. The reason is that, as we showed, at the maximal point once
we have $\nabla_g H=0$ and $\nabla^2_g H$ is negative semi-definite, then the rest of computations in \cite{gv} is the same regardless of the point being in the interior or on the boundary.
\end{pf}

\

As we proved before, there exist two constants
$\bar\delta$ and $C_1$ depending only on $\Vert\nabla Rm_g\Vert$ such that all solutions of
(\ref{eq3}) for some $\delta\leq t\leq t_0$, satisfying
$u_t\leq\bar{\delta}$ satisfy $\parallel\nabla_g
u\parallel_{\infty}<C_1$. Consider now the following quantity:
$$I(M, \partial M, g) :=
\inf_{g'=e^{-2\varphi}g \, , \, \,  \vert \nabla_g \varphi \vert \le C_1\,,\,\,H_{g'}=0}\left(
\int_M R_{g'}^2e^{-\varphi}dV_{g'}\right).$$
We let, for $g'=e^{-2\varphi}g$
$$i(g'):=\int_M R_{g'}^2e^{-\varphi}dV_{g'}.$$
As one can easily check, if two metrics
$g_1$ and $g_2$ are homothetic, then $i(g_1)=i(g_2)$. So, we have
$$I(M, \partial M, g) =
\inf_{g'=e^{-2\varphi}g \, , \, \,  Vol_g'(M)=1 \,\, \hbox {\small and} \,\, \vert \nabla_g \varphi \vert_g \le C_1\,,\,\,H_{g'}=0}\left(
\int_M R_{g'}^2e^{-\varphi}dV_{g'}\right).$$
 Concerning $I(M, \partial M, g)$, we have the following Lemma.
\begin{lem}\label{techn}
There exists a positive constant $C=C(\Vert\nabla Rm_g\Vert)$ such that
$$I(M, \partial M, g) \ge C \left( Y(M,\bM,[g])\right)^2.$$
\end{lem}

\begin{pf} As we have seen
$$I(M,\partial M, g) =
\inf_{g'=e^{-2\varphi}g \, , \, \,  Vol_g' (M)=1 \,\, \hbox {\small and} \,\, \vert \nabla_g \varphi \vert_g \le C_1\,,\,\,H_{g'}=0}\left(
\int_M R_{g'}^2e^{-\varphi}\,dV_{g'}\right).$$
Take $\varphi\in C^\infty(M)$ such that, for $g'=e^{-2\varphi}g$, $Vol_g'(M)=1$ and such that
$\vert \nabla_g \varphi \vert_g \le C_1$ where $C_1$ is given by Proposition \ref{GE}.
Since $Vol_g'(M)=1$, if $p$ is a point where $\varphi$ attains its minimum we have
$$e^{-3\varphi(p)} Vol_g(M) \ge 1,$$
and then, there exists $C_0$ depending only on $(M,g)$ such that
$\varphi(p) \le C_0$. Now, using the mean value theorem, it follows since
$\vert \nabla_g \varphi \vert_g$ is controlled by a constant depending only on
$(M,g)$, that $\max \varphi \le C_0'$ where $C'_0$ depends only on $(M,g)$. Using this, we clearly have that
$$\int_M R_{g'}^2e^{-\varphi}\,dV_{g'} \ge e^{-C_0'} \int_M R_{g'}^2\,dV_{g'}.$$
Using H\"older inequality and the definition of the Yamabe invariant, since $H_{g'}=0$, we get (recall that
$Vol_g' (M)=1$)
$$\int_M R_{g'}^2\,dV_{g'} \ge \left( Y(M,\bM,[g])\right)^2,$$
and then $I(M, \partial M, g) \ge e^{-C_0'}\left( Y(M,\bM,[g])\right)^2$. This ends the proof.
\end{pf}

\

We will prove a lower bound for a solution to the equation (\ref{eq3}) following section 3 in \cite{cd}. Since we are dealing with manifolds with boundary we have to compute the conformal deformation of the integral of $\sigma_2$ in this context. Here is the formula

\begin{lem}\label{change} For a conformal metric $\gt=e^{-2u}g$, we have the following integral transformation
\bea
\int_M \sdtu e^{-4u}\,dV_g&=&\int_M \sdu\,dV_g+\f{1}{8}\int_M R_g|\nabla_g
u|_g^2\,dV_g-\f{1}{4}\int_M|\nabla_g u|_g^4\,dV_g\\
&&+\f{1}{2}\int_M\Delta_g
u|\nabla_g u|_g^2\,dV_g-\f{1}{2}\int_M A_g^1(\nabla_g u,\nabla_g u)\,dV_g\\
&&+\f{1}{4}\oint_{\bM}\dn u \left(R_g+2\Delta_g u-2|\nabla_g u|^2_g\right)\,dS_g\\
&&-\oint_{\bM}A^1_g(\nu,\nabla_g u)\,dS_g-\f{1}{4}\oint_{\bM}\dn |\nabla_g u|^2_g\,dS_g.
\eea
In particular, if the boundary of $M$ is totally geodesic and $\dn u=0$, we get
\bea
\int_M \sdtu e^{-4u}\,dV_g&=&\int_M \sdu\,dV_g+\f{1}{8}\int_M R_g|\nabla_g
u|_g^2\,dV_g-\f{1}{4}\int_M|\nabla_g u|_g^4\,dV_g\\
&&+\f{1}{2}\int_M\Delta_g
u|\nabla_g u|_g^2\,dV_g-\f{1}{2}\int_M A_g^1(\nabla_g u,\nabla_g u)\,dV_g.
\eea
\end{lem}
\begin{pf} For the computations, we will follow section 3 in \cite{cd}. The final formula will be the same as in \cite{cd}, but with some extra terms coming from the boundary.

\

Denote $\sigut=\sutu$, $\sigu=\suu$, $\sigdt=\sdtu$,
$\sigd=\sdu$. We have
$$
2\sigdt=\sigut^2-|A^1_{\gt}|^2_{\gt}\,.
$$
By
equation~\eqref{trasfz2}, we have
$$
\sigut e^{-2u}=\sigu+\Delta_g
u-\f{1}{2}|\nabla_g u|_g^2\,,
$$
so
$$
\sigut^2 e^{-4u}=\sigu^2+(\Delta_g u)^2
+\f{1}{4}|\nabla_g u|_g^4+2\sigu\Delta_g u-\Delta_g u|\nabla_g
u|_g^2-\sigu|\nabla_g u|_g^2\,.
$$
After an easy computation, we get
\bea
|A^1_{\gt}|^2_{\gt}\,\,e^{-4u}&=&|A^1_g|_g^2+|\nabla_g^2
u|_g^2+\f{3}{4}|\nabla_g u|_g^4-\sigu|\nabla_g u|_g^2-\Delta_g u|\nabla_g u|_g^2+\\
&&+\,2(A^1_g)_{ij}\nabla_g^{2\, {ij}} u+2(A^1_g)_{ij}\nabla_g^i u\nabla_g^j
u+2\nabla_{g \, ij}^2 u\nabla_g^i u\nabla_g^j u\,.
\eea
Putting all together,
we obtain
\bea
2\sigdt e^{-4u}&=&2\sigd+(\Delta_g u)^2-|\nabla_g^2
u|_g^2-\f{1}{2}|\nabla_g u|_g^4+2\sigu\Delta_g u\\
&&-\,2(A^1_g)_{ij}\nabla_g^{2\, ij}u-2(A^1_g)_{ij}\nabla_g^i u\nabla_g^j
u-2\nabla^2_{g\, ij} u\nabla_g^i u\nabla_g^j u\,.
\eea
Now, by simple computation, we have the following identities
$$
-2\int_M
(A^1_g)_{ij}\nabla_g^{2\, ij}u\,dV_g=-2\int_M\sigu\Delta_g u\,dV_g+2\oint_{\bM}\dn u\, \sigma_{1}\,dS_{g}-2\oint_{\bM}A^{1}_{g}(\nu,\nabla_{g} u)\,dS_{g}\,,
$$
$$
-2\int_M\nabla^2_{ij} u\nabla_g^i u\nabla_g^j u\,dV_g=\int_M\Delta_g
u|\nabla_g u|_g^2\,dV_g-\oint_{\bM}\dn u\,|\nabla_{g}u|^{2}_{g}\,dS_{g}\,,
$$
where we integrated by parts and we used the Schur's lemma,
$$
2\nabla_g^j (Ric_g)_{ij}= \nabla_i R_g\,,
$$
for the first identity. Finally we get
\bea
2\int_M\sigdt
e^{-4u}\,dV_g&=&2\int_M\sigd\,dV_g\\
&&+\int_M\left[(\Delta_g u)^2-|\nabla_g^2 u|_g^2
-\f{1}{2}|\nabla_g u|_g^4 + \Delta_g u|\nabla_g u|_g^2-2A^1_g(\nabla_g u,\nabla_g
u)\right]\,dV_g\\
&&+\oint_{\bM}\dn u\,\left(\f{1}{2}R_{g}-2A^{1}_{g}(\nu,\nabla_{g} u)-|\nabla_{g}u|^{2}_{g}\right)\,dS_{g}\,,
\eea
Now, integrating the Bochner formula
$$
\f{1}{2}\Delta_{g}|\nabla_{g} u|^{2}_{g}=|\nabla_g^2 u|_g^2+Ric_g(\nabla_g u,\nabla_g
u)\,dV_g+\nabla_{i} u	,\nabla^{i}(\Delta_g u)\,,
$$
we get
$$
\f{1}{2}\oint_{\bM}\dn |\nabla_{g} u|^{2}_{g}\,dS_{g}=\int_{M}\left[|\nabla_{g}^{2}u|^{2}_{g}-(\Delta_{g}u)^{2}+Ric_{g}(\nabla_{g}u,\nabla_{g}u)\right]\,dV_{g}+\oint_{\bM}\dn u\,\Delta_{g} u\,dS_{g}
$$

Using the definition of the Schouten tensor $A^1_g$, we get the first point of the lemma.

\

Now, if the boundary is totally geodesic and $\dn u=0$ on $\bM$, then by Lemma \ref{eq:boun} we have that all the boundary terms must vanish. Thus the second point of the lemma is proved. This completes the proof.
\end{pf}

\

Since $(M,g)$ has totally geodesic boundary, the boundary terms don't effect the conformal transformation of the integral of $\sigma_2$. Hence, following section 3 in \cite{cd} and using Lemma \ref{techn}, we obtain the lower bound.

\begin{pro}[Lower Bound]\label{LB3}
Assume that for some
$t\in[\delta, t_0]$, $t_{0}\leq 2/3$, the following estimate holds
\be
\int_M\sigd(g^{-1}A^1_g)\,dV_g+C\left(\f{7}{10}-t\right)(Y(M,\bM,[g])^2=\mu_t>0,
\ee
for some $C$ depending only on $\Vert\nabla Rm_g\Vert$. Then there exists $\underline{\delta}$ depending only on $diam_g(M)$ and $\Vert\nabla Rm_g\Vert$ such that if $u_t\in
C^2(M)$ is a solution of (\ref{eq3}) and if $A^t_{u_t}\in\Gamma^+_2$
then $u_t\ge\underline{\delta}$.
\end{pro}

\

We have the following  $C^{2,\alpha}$ estimate for solutions of the equation (\ref{eq3}).

\begin{pro}[$C^{2,\alpha}$ estimate]\label{C2}
Let $u_t\in C^4(M)$ be a
solution of (\ref{eq3}) for some $\delta\leq t\leq t_0$, $t_{0}\leq 2/3$, satisfying
$\underline{\delta}<u_t<\bar\delta$, and $\parallel\nabla
u_t\parallel_{{g,\infty}}<C_1$. Then, if $A^t_{u_t}\in\Gamma^+_2$, for $0<\alpha<1$, $\parallel
u_t\parallel_{C^{2,\alpha}}\leq C_2$, where $C_2$ depends only on
$\underline{\delta}, \bar\delta, C_1$ and $\Vert\nabla^2 Rm_g\Vert$.
\end{pro}
\begin{pf}
The interior $C^2$ estimate follows from the work of Chen \cite{ch1} and the boundary $C^2$ estimate follows from Theorem 6 (b) in \cite{ch2}. With the $C^2$ estimate at hand, we obtain high-order estimate (in particular $C^{2,\alpha}$ one) from the works of Evans \cite{ev}, Krylov \cite{kr} and Lions-Trudinger \cite{LT}.
\end{pf}
\

Since we proved $C^{2, \a}$ estimates for solutions of the equation (\ref{eq3}), by the classical Ascoli-Arzela's Theorem, we have that $\mathcal{S}$ is closed, therefore
$\mathcal{S}=[\delta, t_0]$. In particular $t_0\in\mathcal{S}$. Hence the metric \\ $\gt=e^{-2u_{t_0}}g$ then
satisfies $\sigma_2(A^{t_0}_{\gt})>0$, $R_{\gt}>0$ and $L_{\tilde g}=0$. Furthermore, by Lemma \ref{GVineq} we have that the metric $\gt$ satisfies
\be\label{inequa}
(3t_0-2)R_{\gt}\gt<6Ric_{\gt}<3(2-t_0)R_{\gt}\gt.
\ee
Hence the proof of Theorem $\ref{T3_1}$ is complete.

\

Now we are going to give the proof of Theorem $\ref{T3_gen}$.

\

\begin{pfn} {of Theorem $\ref{T3_gen}$}\\
First of all from $R_g>0$ and $L_g=0$, we infer $Y(M,\partial M, [g])>0$. On the other hand, one can easily check that
$$
\sigma_2(g^{-1}A_g)=\frac{3}{16}|R_g|^2-\frac{1}{2}|Ric_g|^2.
$$
Thus, we have $\int_{M}\sigma_2(g^{-1}A_g)\geq 0$ is equivalent to $\int_{M}|Ric_g|^2dV_g\leq \frac{3}{8}\int_{M}|R_g|^2dV_g$. Hence we can apply Theorem $\ref{T3_1}$ with $t_0=\frac{2}{3}$ and get the existence of a metric $\tilde g$ conformal to $g$ such that $Ric_{\tilde g}>0$ and $L_{\tilde g}=0$. Hence appealing to Theorem $\ref{eq:shen}$, we have the proof of Theorem $\ref{T3_gen}$ is complete.
\end{pfn}

\

\section{Four manifolds with boundary}
In this section, we give the proof of Theorem $\ref{T4_gen}$. As for the case of $3$-manifolds, we are going to prove a more general theorem from which Theorem $\ref{T4_gen}$ becomes a direct application.
\begin{thm}\label{T4_1} Let $(M,g)$ be a compact four--dimensional Riemannian
manifold with umbilic boundary and \;$0\leq \alpha\leq 1$. If \;$Y(M,\bM,[g])>0$, and
$$
\frac{1}{2}\kappa_{(P^4,P^3)}-\f{\alpha}{16}\int_M|W_g|^2_g\,dV_g+\f{1}{24}(1-t_0)(2-t_0)Y(M, \partial M, [g])^2>0\,,
$$
for some $t_0\leq 1$, then there exists a conformal metric
$\gt=e^{-2u}g$ whose curvature satisfies
$$
R_{\tilde g}>0,\;\;\;\sigd(\gt^{-1}A^{t_0}_{\gt})-\f{\alpha}{16}|W_{\gt}|^2_{\gt}>0,\hspace{0.1cm}\,\mbox{and}\,\,\,H_{\gt}=0\,.
$$
This implies the pointwise inequalities
$$
(t_0-1)R_{\tilde g}\tilde g< 2Ric_{\tilde g}<(2-t_0)R_{\tilde g}\tilde g.
$$
\end{thm}

\

Throughout the sequel, $(M,g)$ will be a compact
4-dimensional Riemannian manifold with umbilic boundary and with positive Yamabe invariant $Y(M,\bM, [g])$. Since all the hypothesis on the metric $g$ are conformally invariant, then by a result of Escobar, see \cite{es}, we can choose in the conformal class the Yamabe metric, i.e. a metric with positive constant scalar curvature and zero mean curvature. Moreover, since umbilicity is also conformally invariant, we have that the boundary must be totally geodesic. Hence, from now on, $(M,g)$ will be a compact four--manifold with totally geodesic boundary, positive constant scalar curvature and satisfying the integral pinching condition.\\
On the other hand, since $M$ is compact and $R_g>0$,
there exist $t_0>\delta>-\infty$, $\delta<0$ such that $A^\delta_g$ is positive
definite (i.e. $Ric-\f{\delta}{6}R>0$ on $M$). Moreover we can choose $\delta$ so small such that
$$
\sigma_2^{1/2}(g^{-1}A^{\delta}_g)-\f{\sqrt{\alpha}}{4}|W_g|_g>0\,.
$$
Note that $\delta$ depends only on $\Vert Rm\Vert$. \\
Now we define a subclass of the positive cone of order $2$ which will be useful in our arguments

\

\begin{df} For a conformal metric $\gt=e^{-2u}g$, we define the set
$$
\Lambda_{\gt}^+=\left\{t\in[\delta,t_0]\mid
A^t_{\gt}\in\Gamma_2^+\,\,\,\,\mbox{and}\,\,\,\,
\sigma_2^{1/2}(g^{-1}A^t_{\gt})-\f{\sqrt{\alpha}}{4}|W_g|_g>0\right\}\,.
$$
In particular if $t\in \Lambda_{\gt}^+$ then $A^{t}_{\gt}\in\Gamma_{2}^{+}$.
\end{df}

We point out that $\delta\in\Lambda_g^+$.\\\\
\

For $t\in[\delta, t_0]$, consider the path of
equations (in the sequel we use the notation
$A^t_{u_t} := A^t_{g_t}$ for $g_t$ given by
$g_t=e^{-2u_t}g$)
\be\label{eq4}
\left\{
\begin{array}{ll}
\sigma^{1/2}_2(g^{-1}A^t_{u_t})-\f{\sqrt{\alpha}}{4}|W_g|_g=fe^{2u_t}\hspace{1cm} & \mbox{in $M$},\\
\dn u=0\hspace{1cm} &\mbox{on $\bM$}.
\end{array}\right.
\ee
where
$f(x)=\sigma_2^{1/2}(g^{-1}A^{\delta}_g)-\f{\sqrt{\alpha}}{4}|W_g|_g>0$.
Note that $u\equiv 0$ is a solution of (\ref{eq4}) for $t=\delta$.

\

As for the tree-dimensional case, we use the continuity method. Define
$$\mathcal{S}=\left\{t\in[\delta, t_0]\mid\exists\mbox{\,a
solution\,}u_t\in
C^{2,\alpha}(M)\mbox{\,of\,}(\ref{eq4})\mbox{\,with\,} t\in\Lambda^+_{u_t}\right\}.$$
Clearly, with our choice of $f$, $u\equiv 0$ is a solution for
$t=\delta$. Since $\delta\in \Lambda^+_g$, then $\delta\in\mathcal{S}$. Hence, we have $\mathcal{S}\neq\emptyset.$ Let
$t\in\mathcal{S}$, and $u_t$ be a solution. By Proposition
\ref{Inver}, the linearized operator at $u_t$,
$\mathcal{L}^t:C^{2,\alpha}(M)\cap\{\dn u=0\,\,\mbox{on}\,\,\bM\}\rightarrow C^{\alpha}(M)$, is
invertible (note that the additional term in the right hand side of the equation does not effect linearization). The implicit function theorem  tells us
that $\mathcal{S}$ is open. To prove that $\mathcal{S}$ is close we need to establish a priori $C^{2,\a}$ estimates for solutions of the equation (\ref{eq4}). To do so, we start by establishing upper-bound estimate as for the case of $3$-manifolds.

\

\begin{pro}[Upper bound]\label{UB} Let $u_t\in C^2(M)$ be a solution of
(\ref{eq4}) for some $t\in[\delta, t_0]$, with $t\in\Lambda_{u_t}^+$. Then $u_t\leq\bar{\delta}$,
where $\bar{\delta}$ depends only on $\Vert Rm_g\Vert$.
\end{pro}
\begin{pf} From Lemma \ref{PropGamma} (iv), we have
$\f{4}{\sqrt{6}}\sigma_2^{1/2}\leq\sigma_1$, so for all $p\in M$
$$
\f{4}{\sqrt{6}}\f{\sqrt{\alpha}}{4}|W_g|_g+\f{4}{\sqrt{6}}fe^{2u_t}\leq\sigma_1(g^{-1}A^t_{u_t})\,.
$$
Let $p\in M$ be the maximum of $u_t$, then (this is true also if $p\in\bM$, since $\dn u=0$ on $\bM$) we have
$(\Delta u_t)(p)\leq 0$. Then, using (\ref{trasfz}), we have
\bea
\f{4}{\sqrt{6}}\f{\sqrt{\alpha}}{4}(|W_g|_g)(p)+\f{4}{\sqrt{6}}f(p)
e^{2u_t(p)}&\leq&\sigma_1(g^{-1}A^t_{u_t})(p)\\
&=&\sigma_1(g^{-1}A^t_g)(p)+(3-2t)(\Delta u_t)(p)\\
&\leq&\sigma_1(g^{-1}A^t_g)(p)\\
&\leq&\sigma_1(g^{-1}A^\delta_g)(p)\,.
\eea
This implies
$$
\f{4}{\sqrt{6}}f(p)e^{2u_t(p)}\leq\sigma_1(g^{-1}A^\delta_g)(p)
-\f{4}{\sqrt{6}}\f{\sqrt{\alpha}}{4}(|W_g|_g)(p)\,,
$$
where the last term has positive sign.
Since $M$ is compact, this implies $u_t\leq\overline{\delta}$, for some
$\overline{\delta}$ depending only on $\Vert Rm\Vert$.
\end{pf}

\

Following the previous section, once we have an upper bound of the solution, from Proposition \ref{GE}, we get gradient estimates.  Now we are going to establish the lower-bound estimates. To do that we need the following Lemma.
\begin{lem}\label{eq:sigmakap}
 If $\hat g$ is a  Riemannian metric on $M$ conformal to $g$ such that $L_{\hat g}=0$, then
$$
\int_{M}\sigma_2(\hat g^{-1}A_{\hat g})=\frac{1}{2}\kappa_{(P^4,P^3)}.
$$
\end{lem}
\begin{pf}
First of all, one can easily check that the following holds
$$
Q_{\hat g}=2\sigma_2(\hat g^{-1}A_{\hat g})-\frac{1}{12}\D_{\hat g}R_{\hat g} .
$$
Thus integrating this equation and using the divergence theorem, we get
$$
\int_{M} Q_{\hat g}dV_{\hat g}=2\int_{M} \sigma_2(\hat g^{-1}A_{\hat g})dV_{\hat g}+\frac{1}{12}\oint_{\partial M}\frac{\partial R_{\hat g}}{\partial n_{\hat g}}dS_{\hat g}.
$$
On the other hand, since $L_{\hat g}=0$, then
$$
T_{\hat g}=-\frac{1}{12}\frac{\partial R_{\hat g}}{\partial n_{\hat g}}.
$$
Thus we obtain
$$
\int_{M} Q_{\hat g}dV_{\hat g}=2\int_{M} \sigma_2(\hat g^{-1}A_{\hat g})dV_{\hat g}-\oint_{\partial M}T_{\hat g}dS_{\hat g}
$$
Hence, we get
$$
\int_{M}\sigma_2(\hat g^{-1}A_{\hat g})=\frac{1}{2}\kappa_{(P^4,P^3)}.
$$
This completes the proof of the Lemma.
\end{pf}
\

\begin{pro}[Lower bound]\label{LB4} Assume that for some
$t\in[\delta,t_0]$ the following estimate holds
\be
\frac{1}{2}\kappa_{(P^4,P^3)}-\f{\alpha}{16}\int_M|W_g|^2_g\,dV_g+\f{1}{24}(1-t)(2-t)Y(M,\bM,[g])^2=\mu_t>0\,.
\ee
Then there exist $\underline{\delta}$ depending only on
$\mbox{diam}(M,g)$ and $\Vert \nabla^2 Rm\Vert$
such that if $u_t\in C^2(M)$ is a solution of~\eqref{eq4} and if
$t\in\Lambda_{u_t}^+$ then $u_t\geq\underline{\delta}$.
\end{pro}
\begin{pf} Since $A^t_g=A^1_g+\f{1-t}{2}\sigma_1(A^1_g)g,$ we
easily have
$$
\sigma_2(A^t_g)=\sigma_2(A^1_g)+\f{3}{2}(1-t)(2-t)\sigma_1(A^1_g)^2\,.
$$
Letting $\gt=e^{-2u_t}g$, since $u_{t}$ is a solution of equation~\eqref{eq4}, we have
$$
f^2e^{4u_t}+\f{\sqrt{\alpha}}{2}f|W_g|_ge^{2u_t}=\sigma_2(g^{-1}A^t_{u_t})-\f{\alpha}{16}|W_g|_g^2\,.
$$
The left--hand side can be estimate by
$$
f^2e^{4u_t}+\f{\sqrt{\alpha}}{2}f|W_g|_ge^{2u_t}\leq C'e^{2u_t}\,,
$$
where the positive constant $C'$ depends only on
$\Vert Rm\Vert$. So we get
\bea
C' e^{2u_t}&\geq&\sigma_2(g^{-1}A^t_{u_t})-\f{\alpha}{16}|W_g|_g^2\\
&=&e^{-4u_t}\left(\sigma_2(\gt^{-1}A^1_{u_t})+\f{1}{24}(1-t)(2-t)R_{\gt}^2\right)-\f{\alpha}{16}|W_g|_g^2\,.
\eea
Integrating this with respect to $dV_g$, we obtain
\bea
C'\int_Me^{2u_t}\,dV_g&\geq&\int_M\sigma_2(\gt^{-1}A^1_{u_t})\,dV_{\gt}
-\f{\alpha}{16}\int_M|W_g|_g^2\,dV_g+\f{1}{24}(1-t)(2-t)\int_MR_{\gt}^2\,dV_{\gt}\\
&=&\frac{1}{2}\kappa_{(P^4,P^3)}
-\f{\alpha}{16}\int_M|W_g|_g^2\,dV_g+\f{1}{24}(1-t)(2-t)\int_MR_{\gt}^2\,dV_{\gt}\\
&\geq&\frac{1}{2}\kappa_{(P^4,P^3)}-\f{\alpha}{16}\int_M|W_g|_g^2\,dV_g
+\f{1}{24}(1-t)(2-t)Y(M, \partial M, [g])^2=\mu_t>0\,,
\eea
where we have  used  Lemma $\ref{eq:sigmakap}$, and the fact that  for any conformal metric $g'\in[g]$, if $H_{g'}=0$, then
$$
\int_MR^2_{g'}\,dV_{g'}\geq Y(M,\bM,[g])^2\,.
$$
This gives
$$
\max_M u_t\geq\log{\mu_t}-C(\mbox{diam}(M,g),\Vert Rm\Vert)\,.
$$
Since, as already remarked $\max_{M}| \nabla_g u_t |_{g} \leq C_1$ by the same arguments as the ones of Proposition~\ref{GE} , then we have the Harnack inequality
$$\max_M u_t\leq\min_M u_t + C(\mbox{diam}(M,g),\Vert \nabla^2 Rm\Vert)\,,
$$ by simply integrating along
a geodesic connecting points at witch $u_t$ attains its maximum
and minimum. Combining this two inequalities, we obtain
$$
u_{t}\geq\min_M
u_t \geq\log{\mu_t}-C=:\underline{\d}\,,
$$
where $C$ depends only on $\mbox{diam}(M,g)$ and $\Vert \nabla^2
Rm\Vert$.
\end{pf}

Once we have $C^{0}$ and $C^{1}$ estimates, using the same arguments as the ones of Proposition \ref{C2}, we get $C^{2,\a}$ estimates. Thus we are ready to apply the continuity method as in the $3$-dimensional case, and conclude the proof of Theorem~\ref{T4_1}.

\

Now we are ready to present the proof of Theorem $\ref{T4_gen}$.

\

\begin{pfn} {of Theorem $\ref{T4_gen}$}\\
First of all, since $Y(M, \partial M, [g])>0$, and $\kappa_{(P^4,P^3)}>\frac {1}{8}\int_{M}|W_g|^2dV_g$, then we can apply Theorem $\ref{T4_1}$ with $t_0=1$ and $\alpha=1$ and get the existence of a metric $\tilde g$ conformal to $g$ such that
$$
\sigma_{2}(g^{-1}A_{\tilde g})>\frac{1}{16}|W_g|^2, \;\;\text{and}\;\; L_{\tilde g}=0.
$$
This is equivalent to
$$
\sigma_{2}(\tilde g^{-1}A_{\tilde g})>\frac{1}{16}|W_{\tilde g}|^2, \;\;\text{and } L_{\tilde g}=0.
$$
On the other hand, one can check easily that the following holds
$$
\sigma_{2}(\tilde g^{-1}A_{\tilde g})=\frac{1}{96}R^2_{\tilde g}-\frac{1}{8}|E_{\tilde g}|^2.
$$
Thus, we obtain
$$
\frac{1}{6}R^2_{\tilde g}-2|E_{\tilde g}|^2>|W_{\tilde g}|^2
$$
So rearranging the latter inequality, we get the Margerin's weak pinching condition, namely
$$
WP_{\tilde g}<\frac{1}{6}.
$$
Hence, applying Theorem $\ref{eq:bmar}$, we conclude the proof of Theorem $\ref{T4_gen}$.
\end{pfn}

\

\section{Principal eigenvalue of $P^{4,3}_g$ and applications}
In this section, we provide the proof of Theorem $\ref{eq:spect}$ and Corollary $\ref{eq:existence}$. We start by giving a Proposition which will be used for the proof of Theorem $\ref{eq:spect}$.
\begin{pro}\label{eq:pk}
Let $(M,g)$ be a compact four--dimensional Riemannian manifold with boundary such that $L_g=0$. Assuming $Ric_g\leq R_gg$, then we have $P^{4,3}_g$ is a non-negative operator and $\ker P^{4,3}_g \simeq R$.
\end{pro}
\begin{pf}
First of all, since $L_g=0$, then for every $u\in H_{\frac{\partial }{\partial n}}$, we have
$$
\left<P^{4,3}_gu,u\right>_{L^2(M)}=\int_{M}(\D_g u)^2dV_g+\frac{2}{3}\int_MR_g|\n_gu|^2dV_g-2\int_MRic_g(\n_gu,\n_gu)dV_g.
$$
Now we recall  the Bochner identity
$$
\frac{1}{2}\D_g(|\n_gu|^2)=|\n^2_gu|^2+Ric_g(\n_gu,\n_gu)+<\n_gu,\n_g(\D_gu)>.
$$
Integrating the latter formula, applying the divergence theorem and integration by part, we get
$$
-\frac{1}{2}\oint_{\partial M}\frac{\partial (|\n_gu|^2)}{\partial n_g}dS_g=\int_{M}|\n^2_gu|^2dV_g+\int_MRic_g(\n_gu,\n_gu)-\int_M(\D_gu)^2dV_g-\oint_{\partial M}\frac{\partial u}{\partial n_g}\D_gudS_g.
$$
Recalling that $u\in H_{\frac{\partial }{\partial n}}$, then
$$
\oint_{\partial M}\frac{\partial u}{\partial n_g}\D_gudS_g=0.
$$
Using Lemma $\ref{eq:boun}$, we obtain
$$
\frac{\partial (|\n_gu|^2)}{\partial n_g}=0.
$$
Hence, we get
$$
\int_{M}|\n^2_gu|^2dV_g+\int_MRic_g(\n_gu,\n_gu)=\int_M(\D_gu)^2dV_g.
$$
Now, using the latter formula, we have
$$
\left<P^{4,3}_gu,u\right>_{L^2(M)}=-\frac{1}{3}\int_{M}(\D_g u)^2dV_g+\frac{4}{3}\int_{M}|\n_g^2u|^2dV_g+\frac{2}{3}\int_MR_g|\n_gu|^2dV_g-\frac{2}{3}\int_MRic_g(\n_gu,\n_gu)dV_g.
$$
Next, setting
$$
\bar{\n_g}^2u=\n_g^2u-\frac{1}{4}\D_g g;
$$
we get
$$
\left<P^{4,3}_gu,u\right>_{L^2(M)}=\frac{4}{3}\int_{M}|\bar \n_g^2u|^2dV_g+\frac{2}{3}\int_M(R_gg-Ric_g)(\n_gu,\n_gu)dV_g.
$$
So using the hypothesis $R_gg-Ric_g\geq 0$, we infer
$$
\left<P^{4,3}_gu,u\right>_{L^2(M)}\geq\frac{4}{3}\int_{M}|\bar \n_g^2u|^2dV_g.
$$
Hence, we obtain $P^{4,3}_g$ is a non-negative operator. So to finish the proof of the Proposition, it remains only to show that the kernel is constituted only by constants. In order to do that, we assume that there exists a non constant function $u\in H_{\frac{\partial}{\partial n}}$ such that $P^{4,3}_gu=0$, and argue for a contradiction.  From the fact that  $u\in ker P^{4,3}_g$, we infer that
$$
\n_g^2u-\frac{1}{4}\D_g g=0.
$$
Now calling the doubling of $M$ by $DM$, and the reflected metric by $\bar g$,  we have that $\bar g$ is $C^{2,\alpha}$. Next  we reflect $u$ across $\partial M$ and call the reflection by $u_{DM}$. Thus, we obtain an element in $H^2(DM)$ verifying
$$
\n_{\bar g}^2u_{DM}-\frac{1}{4}\D_{\bar g} \bar g=0.
$$
Thus using a result of Tashiro \cite{Ta}, we infer that $(DM, \bar g)$ is conformally diffeomorphic to $S^4$. Thus $(M,g)$ is also conformally diffeomorphic to $S^4_+$. So we derive  the existence of a metric $\tilde g$ conformal to $g$ on $M$ which is Einstein, of constant positive scalar curvature, and $L_{\tilde g}=0$. Hence using the conformal invariance of $P^{4,3}_g$ , we get
$$
\frac{4}{3}\int_{M}|\bar \n_{\tilde g}^2u|^2dV_{\tilde g}+\frac{1}{2}R_{\tilde g}\int_{M}|\n_{\tilde g} u|^2dV_{\tilde g}=0.
$$
Thus, we obtain $u$ is constant and reach a contradiction. This completes the proof of the Proposition.
\end{pf}

\

Having this at hand, we are ready to give the proof of Theorem $\ref{eq:spect}$.

\

\begin{pfn}{ of Theorem $\ref{eq:spect}$}\\
Applying Theorem $\ref{T4_1}$ with $t_0=0$, and $\alpha=0$, we get the existence of a metric $\tilde g$ conformal to $g$ such that
$$
Ric_{\tilde g}\leq R_{\tilde g}\tilde g,\;\;\;\text{and}\;\;\;L_{\tilde g}=0.
$$
Hence appealing to Proposition $\ref{eq:pk}$, we obtain that $P^{4,3}_{\tilde g}$ is non-negative and $\ker P^{4,3}_{\tilde g}\simeq \R$. Now recalling that the non-negativity of the operator $P^{4,3}_g$ and the triviality of its kernel are conformally invariant properties, we have that the proof of Theorem $\ref{eq:spect}$ is complete.
\end{pfn}

\

Next, we are going to present the proof of corollary $\ref{eq:existence}$.

\

\begin{pfn} { of Corollary $\ref{eq:existence}$}\\
Due to $\eqref{eq:lawqt}$, the existence of constant $Q$-curvature, constant $T$-curvature ad zero mean curvature metrics conformal to the background one $g$ is equivalent to solving the following (BVP)
\begin{equation}\label{eq:bvps}
\left\{
\begin{split}
P^4_gu+2Q_g&=2\bar Qe^{4u}\;\;&\text{in}\;\;M;\\
P^3_gu+T_g&=\bar Te^{3u}\;\;&\text{on}\;\;\partial M;\\
\frac{\partial u}{\partial n_g}&=0\;\;&\text{on}\;\;\partial M.
\end{split}
\right.
\end{equation}
where $\bar Q$ and $\bar T$ are constant real numbers. On the other hand it is easy to see that critical points of the functional
$$
II(u)=\left<P^{4,3}_g u, u\right>_{L^2(M)}+4\int_MQ_g udV_g+\oint_{\partial M}T_g udS_g-\kappa_{P^4_g}\log\int_Me^{4u}dV_g-\frac{4}{3}\left(\kappa_{(P^4,P^3)}-\kappa_{P^4_g}\right)\log\oint_{\partial M}e^{3u}dS_g;
$$
are weak solution of $\eqref{eq:bvps}$, hence from standard elliptic regularity theory, they are smooth solutions. Thus to prove the corollary, we will prove the existence of critical points. More precisely, under our assumption, we will prove the existence of a minimizer. To do so, we first point out that the functional $II$ is invariant by translation by constant, and can also be written in the following form
\begin{equation}
\begin{split}
II(u)=\left<P^{4,3}_g u, u\right>_{L^2(M)}+4\int_M Q_g (u-\bar u)dV_g+\oint_{\partial M}T_g(u-\bar u_{\partial M})dS_g-\kappa_{P^4_g}\log\int_Me^{4(u-\bar u)}dV_g\\-\frac{4}{3}\left(\kappa_{(P^4,P^3)}-\kappa_{P^4_g}\right)\log\oint_{\partial M}e^{(3u-\bar u_{\partial M})}dS_g;
\end{split}
\end{equation}
Now exploiting this way of writing $II$, we have if $\kappa_{P^4_g}\leq 0$, then by Jensen's inequality, we obtain
$$
II(u)\geq \left<P^{4,3}_g u, u\right>_{L^2(M)}+4\int_MQ_g (u-\bar u)dV_g+\oint_{\partial M}T_g (u-\bar u_{\partial M})dS_g
$$
Hence, using Cauchy inequality, trace theorem, Sobolev embedding, Poincar\'e inequality, and Lemma $\ref{eq:neq}$, we get
$$
II(u)\geq\gamma ||u-\bar u||_{H^2}-C
$$
for some $\gamma>0$ and some large $C$. Next if $\kappa_{P^4_g}>0$, we use Proposition $\ref{eq:mos-tru}$ and Proposition $\ref{eq:mos-trub}$ to obtain
\begin{equation}
\begin{split}
II(u)\geq\left<P^{4,3}_g u, u\right>_{L^2(M)}+4\int_M Q_g (u-\bar u)dV_g+\oint_{\partial M}T_g(u-\bar u_{\partial M})dS_g\\+\left(-\frac{4}{\alpha_1}\kappa_{P^4_g}-\frac{3}{\alpha_2}\left(\kappa_{(P^4,P^3)}-\kappa_{P^4_g}\right)\right)\left<P^{4,3}_g u, u\right>_{L^2(M)}-C_{\alpha_1,\alpha_2};
\end{split}
\end{equation}
for  $\alpha_1<16\pi^2$ and $\alpha_2<12\pi^2$, and $C_{\alpha_,\alpha_2}$ a constant depending only on $\alpha_1$, $\alpha_2$ and $(M,g)$. To continue the proof we need the following rigidity result
\begin{lem}\label{eq:rigid}
Let $(M,g)$ be a compact four-dimensional Riemannian manifold with umbilic  boundary. Assuming that \;$Y(M,\partial M,[g])\geq0$, we have  \;$\kappa_{(P^4,P^3)}\leq 4\pi^2$ and equality holds if and only if \;$(M,g)$ is conformally diffeomorphic to \;$S^4_+$ with its standard metric
\end{lem}
\begin{pf}
Since $Y(M, \partial M, [g])\geq 0$ and $(M,g)$ has umbilic boundary, then by a result of Escobar \cite{es}, we can take the Yamabe metric $\tilde g$ which has constant non-negative scalar curvature and such that $L_{\tilde g}=0$. On the other hand, still by a result of Escobar \cite{es}, we have  that
$$
Y(M,\partial M, [g])=R_{\tilde g}Vol_{\tilde g}(M)^{\frac{1}{2}}\leq Y(S^4_+, S^3, [g_S])=8\sqrt{3}\pi;
$$
and equality holds if and only if $(M,g)$ is conformally diffeomorphic to $(S^4_+,g_{S})$. Now using Lemma $\ref{eq:sigmakap}$, we have
$$
\kappa_{(P^4,P^3)}=2\int_M\sigma_2(\tilde g^{-1}A_{\tilde g})dV_{\tilde g}.
$$
On the other hand, we have also
$$
\sigma_2(\tilde g^{-1}A_{\tilde g})=\frac{1}{96}R^2_{\tilde g}-\frac{1}{8}|E_{\tilde g}|^2.
$$
Thus, we arrive to
$$
\kappa_{(P^4,P^3)}=\frac{1}{4}\int_M\left(\frac{1}{12}R^2_{\tilde g}-|E_{\tilde g}|^2dV_{\tilde g}\right)\leq \frac{1}{48}R_{\tilde g}^2Vol_{\tilde g}(M)\leq \frac{192}{48}\pi^2=4\pi^2
$$
and equality holds if and only if $(M,g)$ is conformally diffeomorphic to $(S^4_+,g_{S})$. This completes the proof of the lemma.
\end{pf}

\

Now coming back to our proof, we have that, since  $Y(M,\partial M,[g])>0$, and $(M,g)$ has an umbilic boundary, then by Lemma $\ref{eq:rigid}$  $\kappa_{(P^4,P^3)}\leq 4\pi^2$ and equality holds if and only if $(M,g)$ is conformally diffeomorphic to $S^4_+$ with its standard metric. Hence, we can assume that $\kappa_{(P^4,P^3)}<4\pi^2$, otherwise there is noting to do. Thus taking $\alpha_1$ close to $16\pi^2$ and $\alpha_2$ close to $12\pi^2$, and using Cauchy inequality, trace theorem, Sobolev embedding, Poincar\'e inequality, and Lemma $\ref{eq:neq}$, we get
$$
II(u)\geq\gamma _0||u-\bar u||_{H^2}-C_0;
$$
for some $\gamma_0>0$ and some large $C_0$. Hence in any case we obtain
\begin{equation}\label{eq:coercive}
II(u)\geq\gamma _||u-\bar u||_{H^2}-C_1
\end{equation}
for some $\gamma_1>0$ and some large $C_1$. From this, and the fact that $II$ is invariant by translation by constant, we have the existence of a minimizer $u_n$ such that
\begin{equation}\label{eq:normal}
\int_Me^{4u_n}dV_g=1.
\end{equation}
 Thus by the coercivity property $\eqref{eq:coercive}$, we have
$$
||u_n-\bar u_n||_{H^2}\leq C.
$$
On the other hand, using  Proposition $\ref{eq:mos-tru}$, we infer
\begin{equation}\label{eq:log}
\log \int_Me^{4(u_n-\bar u_n)}dV_g\leq C
\end{equation}
So using $\eqref{eq:normal}$, $\eqref{eq:log}$ and Jensen's inequality we infer
$$
|\bar u_n|\leq C.
$$
Thus, we arrive to
\begin{equation}\label{eq: bounded}
||u_n||_{H^2}\leq C
\end{equation}
Hence up to a subsequence, we have
$$
u_n\rightharpoonup u \\;\;\text{in}\;\;H^2.
$$
Furthermore, we have $u\in H_{\frac{\partial }{\partial n}}$. On the other hand, it is easy to see that $II$ is weakly lower semicontinuous on $H^2$. Thus we have $u$ is a minimizer of $II$. This completes the proof of Corollary $\ref{eq:existence}$.
\end{pfn}

\


\begin{thebibliography}{99}



\bibitem{bes}  Besse A. L., {\em Einstein Manifolds}, Springer Verlag, Berlin, 1987.


\bibitem{bran1} Branson T., {\em The functional determinant}, Global Analysis Research Center Lecture Note Series, Number 4, Seoul National University (1993).

\bibitem{bran2} Branson T., {\em Differential operators canonically associated to a conformal structure}, Math. Scand. 57 (1995), no. 2, 293-345.

\bibitem{bo} Branson T.P., Oersted B., {\em Explicit functional determinants in four dimensions}, Proc. Amer. Math. Soc. 113 (1991), no. 3, 669-682.

\bibitem{be} Berger M., {\em Les varietes  Riemanniennes $\frac{1}{4}$-pincees}. Ann. Scuola Norm. Sup. Pisa 14 (1960), 161-170.

\bibitem{bs}  Brendle S., Schoen R. M., {\em Classification of manifolds with weakly $\frac{1}{4}$-pinched curvatures},  Acta Math.  200  (2008),  no. 1, 1--13.

\bibitem{bw} Bohm. C., Wilking B., {\em Manifolds with positive curvature operator are space forms}, Ann. of Math. 167 (2008), no. 3, 1079Ð1097.
            
\bibitem{cns} Caffarelli L., Nirenberg L., Spruck J., {\em  The Dirichlet problem for nonlinear second order elliptic equations}. III. Functions of the eigenvalues of the Hessian, Acta Math. 155 (1985), no. 3-4, 261-301.

\bibitem{cd} Catino G., Djadli Z., {\em  Conformal deformations of integral pinched 3-manifolds}, to appear in Adv. Math.

\bibitem{CDN} Catino G., Djadli Z., Ndiaye C.B., {\em A sphere theorem on locally conformally flat even--dimensional manifolds}, preprint, 2008.

\bibitem{ch1} Chen S. S., {\em Local estimates for some fully nonlinear elliptic equations}, Int. Math. Res. Not. 55 (2005), 3403-3425.

\bibitem{ch2} Chen S. S., {\em Conformal deformation on manifolds with boundary}, Preprint Server: arXiv:0811.2521v1, 2008.

\bibitem{cq1} Chang S.Y.A., Qing J.,{\em The Zeta Functional Determinants on manifolds with boundary 1. The Formula}, J. Funct. Anal. 147 (1997), 327-362.

\bibitem{cq2} Chang S.Y.A., Qing J.,{\em The Zeta Functional Determinants on manifolds with boundary II. Extremal Metrics and Compactness of Isospectral Set}, Journal of Functional Analysis 147 (1997), 363-399.

\bibitem{cgy2} Chang S.Y.A., Gursky M.J., Yang P.C.,{\em A conformally invariant sphere theorem in four dimensions}, Publ. Math. Inst. Hautes etudes Sci. 98 (2003), 105-143.

\bibitem{cy1}  Chang S.Y.A., Yang P.C., {\em On a fourth order curvature invariant}, Comtemporary Mathematics, 237, Spectral Problems in Geometry and Arithmetic, Ed. T. Branson, AMS, 1999, 9-28.

\bibitem{dm} Djadli Z., Malchiodi A., {\em Existence of
conformal metrics with constant $Q$-curvature}, Ann. of Math. 168 (2008), 813Ð858.

\bibitem{es} Escobar J.F., {\em The Yamabe problem on manifolds with boundary}, J. Differential Geom. 35 (1992), no.1, 21-84.

\bibitem{ev} Evans L. C., {\em Classical solutions of fully nonlinear, convex, second-order elliptic equations}. Comm. Pure Appl. Math. 35 (1982), no. 3, 333-363.

\bibitem{fg} Fefferman C., Graham C,R., {\em $Q$-curvature and Poincar\'e metrics},  Math. Res. Lett.  9  (2002),  no. 2-3, 139--151.


\bibitem{fg1} Fefferman C., Graham C., {\em Conformal invariants}, In Elie Cartan et les mathematiques d'aujourd'hui. Asterisque (1985), 95-116.


\bibitem{ha} Hamilton R. S., {\em  Three manifolds with positive Ricci curvature}, J. Differential Geom. 17 (1982), no. 2, 255-306.

\bibitem{gt} Gilbar D., Trudinger N., {\em  Elliptic Partial Differential Equations of Second Order}, 2nd edition, Springer-Verlag, 1983.

\bibitem{gjms} Graham C,R., Jenne R., Mason L., Sparling G., {\em Conformally invariant powers of the laplacian, I: existence},  J. London Math. Soc. 46 (1992), no.2, 557-565.

\bibitem{gv} Gursky M. J.,  Viaclovsky J. A., {\em A fully nonlinear equation on four-manifolds with positive scalar curvature},  J. Differential Geom. 63 (2003), no. 1, 131-154.

 \bibitem{kl} Klingenberg W., {\em Uber Riemannsche Mannigfaltigkeiten mit nach oben beschrankter Krummung}. Ann. Mat. Pura Appl., 60 (1962), 49-59.

\bibitem{kr} Krylov N. V., {\em Boundedly inhomogeneous elliptic and parabolic equations in a domain}, Izv. Akad. Nauk SSSR Ser. Mat. 47 (1983), no. 1, 75-108.

\bibitem{LT} Lions P.L., Trudinger N.S., {\em Linear oblique derivative problems for the uniformly elliptic Hamilton-Jacobi-Bellman equation}, Math. Z. 191 (1986), no. 1, 11-15.

\bibitem{ma} Margerin C., {\em  A sharp characterization of the smooth 4-sphere in curvature terms}, Comm. Anal. Geom. 6 (1998), 21-65.

\bibitem{nd1} Ndiaye C.B., {\em Conformal metrics with constant $Q$-curvature for manifolds with boundary},  Comm. Anal. Geom.  16  (2008),  no. 5, 1049--1124.

\bibitem{nd2} Ndiaye C.B., {\em Constant \;$T$-curvature conformal metric on 4-manifolds with boundary},   Pacific J. Math.  240  (2009),  no. 1, 151--184.

\bibitem{p1} Paneitz S., {\em A quartic conformally covariant differential operator for arbitrary pseudo-Riemannian manifolds}, preprint, 1983.

\bibitem{p2} Paneitz S., {\em Essential unitarization of symplectics and applications to field quantization}, J. Funct. Anal. 48 (1982), no. 3, 310-359.

\bibitem{ra}  Rauch H. E., {\em A contribution to differential geometry in the large}, Ann. of Math. 54 (1951), 38-55.

\bibitem{via} Viaclovsky J., {\em Conformal geometry, contact geometry, and the calculus of variations}, Duke Math. J. 101 (2000), no. 2, 283-316.

\bibitem{she} Shen Y., {\em  On Ricci deformation of a Riemannian metric on manifold with boundary}, Pacific J. Math.  173  (1996),  no. 1, 203--221

\bibitem{Ta} Tashiro, Y., {\em Complete Riemannian manifolds and some vector fields}, Trans. Amer. Math. Soc. 117 (1965), 251-275.
\end{thebibliography}
\end{document}